\documentstyle[11pt,amssymb,amsmath,amstext,amsfonts,]{article}

\newcommand{\sq}{{$\enspace\square$}}

\newcommand{\onto}{\,\twoheadrightarrow\,}
\newcommand{\ad}{{\mathtt{ad}^{\,}}}
\newcommand{\rk}{{\mathtt{rk}^{\,}}}

\newcommand{\g}{{\mathfrak g}}
\newcommand{\uu}{{\mathfrak u}}
\newcommand{\p}{{{\mathfrak p}}}

\newcommand{\NN}{{\mathcal{N}\!}\mbox{\it{ilp}}}

\newcommand{\bg}{{\mathfrak b}}

\newcommand{\n}{{\cal N}}
\newcommand{\PP}{{\bold P}}
\newcommand{\pr}{{\mathtt{pr}}}

\newcommand{\Ad}{{\mathtt{Ad}^{\,}}}

\newcommand{\C}{{\Bbb C}}

\newcommand{\inv}{^{-1}}

\newcommand{\ab}{{\hskip 8mm}}

\newcommand{\hh}{\mbox{\it Hitch}}

\newcommand{\X}{{\widetilde M}}

\newcommand{\sset}{\subset}
\newcommand{\Spec}{\mbox{Spec}}

\begin{document}
\setlength{\parindent}{0pt}
\setlength{\parskip}{3pt plus 5pt minus 0pt}

\centerline{\Large {\bf The global nilpotent variety is Lagrangian}}

\vskip 10mm
\centerline{\Large {Victor Ginzburg}}
\vskip 1cm

\ab The purpose of this note is to present a short elementary
proof of a theorem due to Faltings and Laumon,
 saying that the global nilpotent
cone is a Lagrangian substack in the cotangent bundle of the
moduli space of $G$-bundles on a complex compact curve. This result
plays a crucial role in the Geometric Langlands program, see [BD],
since it insures that the ${\mathcal{D}}$-modules on the
moduli space of $G$-bundles whose characteristic variety
is contained in  the global nilpotent
cone are automatically {\it holonomic}, hence, e.g. have
finite length.
\medskip

\ab  Let $(M, \omega)$ be a smooth
symplectic
algebraic variety. A (possibly singular) algebraic subvariety $Y
\subset M$ is said to be
isotropic, resp. Lagrangian, if the tangent space, $T_yY$, at any
regular
point $y\in Y$ is an isotropic, resp. Lagrangian,
 vector subspace in the symplectic
vector space $T_yM$ (we always assume $Y$ to be {\it reduced},
but not necessarily irreducible). The following characterisation
of isotropic subvarieties proved, e.g. in \cite[Prop. 1.3.30]{CG},
will be used later
: {\it $Y\subset M$ is isotropic if and only if for any smooth 
locally closed subvariety $W \subset Y$, we have} $\omega|_{_W}=0\,$
(here $W$ is 
possibly contained in the singular
locus of $Y$). An advantage of this 
characterisation is that it allows to extend the notion of `being
isotropic' from algebraic subvarieties to semi-algebraic constructible
subsets. Thus, we call a semi-algebraic constructible
subset $Y\subset M$ isotropic if $\omega|_{_W}=0$ for any smooth 
locally closed algebraic variety $W \subset Y$.

\ab Now, let $M$ be a smooth stack
that can be
locally presented
 as $p: \X \onto M$, where $\X$ is a smooth algebraic
variety and $p$ is a smooth surjective morphism; e.g., 
$M$ is locally isomorphic to the quotient
of a smooth algebraic variety 
modulo an algebraic action of an algebraic group, see e.g. [LMB].
We have a natural
diagram
\[ T^*M\;\twoheadleftarrow\; T^*M\times_M \X\;
\stackrel{\epsilon}{\hookrightarrow}  \;             T^*\X.\] 
A substack $Z\sset T^*M$ is said to be constructible,
resp. isotropic or Lagrangian,
if $\,\epsilon(Z\times_M \X)\,$ is a constructible,
resp. isotropic or Lagrangian, subset of $T^*\X$ relative
to the standard symplectic structure on the cotangent bundle
of a smooth variety.

\ab Let $X$ be a smooth complex compact connected 
algebraic curve of genus
$g>1$, and $G$ a complex 
semisimple\footnote{it is not hard to extend our results
to any reductive group, but that would lead to
unpleasant dimension shifts in various formulas below,
so we restrict ourselves to the semisimple case.} group. Below, 
we write $Bun_G$ for the moduli space of principal algebraic
$G$-bundles on $X$, regarded as a {\it stack}, cf. \cite{LMB},
rather than a scheme. In particular, no stability conditions
on principal
bundles $P\in Bun_G$ are imposed. 
Given a principal $G$-bundle
$P$, let ${\mathfrak g}_{_P}$ and ${\mathfrak g}^*_{_P}$
denote the associated vector bundles corresponding to the adjoint and
co-adjoint representations of $G$, respectively. 
Let $\PP$ be the universal bundle on $Bun_G\times X$, and
$p: Bun_G\times X\to Bun_G$ the projection.

\ab The
cotangent stack,
$T^*Bun_G,$  is a stack, see [BD, n.1.1.1], which is
relatively representable over $Bun_G$
by the affine spectrum of the sheaf of algebras 
${\mathtt{Sym}} (R^1\!p_*{\mathfrak g}_{_\PP})$,
the Symmetric algebra of the first derived pushforward sheaf.
 Note that for all $i>1$, we have: $R^i\!p_*{\mathfrak g}_{_\PP}=0$,
since $\dim X=1$. Hence the formation $R^1\!p_*(-)$
is right-exact, and therefore commutes with base change.
For a scheme $S$ and a morphism:
$S \to Bun_G$, write $\PP(S)$ for the pull-back of the universal bundle to
$S\times X$. Using the base change, one obtains the following
(stack version of the) Kodaira-Spencer formula for the
set  of $S$-points of the fiber of
$T^*Bun_G$ over $\PP(S)$:
\begin{align*}
T^*_{\PP(S)}Bun_G & =\Gamma\bigl(S,\,
{\mathcal{H}}^{\!}om_{_{{\mathcal{C}^{\!}}oh}}
(R^1\!p_*{\mathfrak g}_{_{\PP(S)}}\,,\,{\mathcal{O}}_{_S})
\bigr)\\
&=\Gamma\bigl(S,\,{\mathcal{H}}^{\!}om_{_{D^b({\mathcal{C}^{\!}}oh)}}
(Rp_*{\mathfrak g}_{_{\PP(S)}}[1]\,,\,{\mathcal{O}}_{_S})
\bigr)\hspace{37mm}(1)\\
&=\Gamma\bigl(S,\,
{\mathcal{H}}^{\!}om_{_{D^b({\mathcal{C}^{\!}}oh)}}({\mathfrak
g}_{_{\PP(S)}}\,,\,Rp^!{\mathcal{O}}_{_S}[-1])
\bigr)\\
&=\Gamma(S\times X,\,\, {\mathfrak g}^*_{_{\PP(S)}}\otimes
\Omega^1_{X\times S/S})\,
\simeq\,\Gamma(S\times X,\, \,{\mathfrak g}_{_{\PP(S)}}\otimes
\Omega^1_{X\times S/S})\,,
\end{align*}
where the second isomorphism exploits the fact that
the complex $Rp_*{\mathfrak g}_{_{\PP(S)}}[1]$ is concentrated
in non-positive degrees, and the last isomorphism uses the identification 
${\mathfrak g}^*_{_\PP}\simeq{\mathfrak g}_{_\PP}$ induced by the Killing form on 
${\mathfrak g}$.

\ab Write $\n$ for the nilpotent cone in ${\mathfrak g}$,
the zero variety of the set of $\Ad G$-invariant polynomials on
 ${\mathfrak g}$ without constant term.
 Choose a Borel subgroup $B\subset G$ with Lie algebra 
${\mathfrak b}$, and let ${\mathfrak n}$ denote the nilradical of ${\mathfrak b}$.
Kostant proved [Ko], see also [CG, ch.6], that  $\n$ is
equal, as a subscheme of $\g$, to the image of the Springer resolution, the
morphism: $G\times_{_B}{\mathfrak n} \to {\mathfrak g}$, given by
the assignment: $(g, x) \mapsto \Ad g(x)$.

\ab Given a scheme $S$ and a $G$-bundle $P$ on $S\times X$, we
choose  local trivialisations of the vector bundles
${\mathfrak g}_{_{P}}$ and $\Omega^1_X$ and view a local section $x$
of ${\mathfrak g}_{_{P}}\otimes
\Omega^1_X$ as a function $S\times X\to {\mathfrak g}$. The section $x$ is called
{\it nilpotent} if the corresponding function 
gives a morphism $S\times X\to \n\sset{\mathfrak g}$.
The notion of a nilpotent section
 does not depend on the choices of trivialisations involved.
 Following Laumon \cite{La}, define the {\it global nilpotent
cone} as a closed (non-reduced) substack
$\NN \,\subset \, T^*Bun_G$ whose set of $S$-points is 
\[\NN(S)=\{({P},x)\;\;\big|\;\; x\in \Gamma(S\times X,\,{\mathfrak g}_{_{P}}\otimes
\Omega^1_X)\,\;, \;\,x\;\mbox{{\small {\it is nilpotent section}}}\},\]
where $P$ runs over $G$-bundles on $S\times X$.
\vspace{2mm}

{\bf Main Theorem.} $\;\NN$ {\it is a Lagrangian substack in} $T^*Bun_G$.
\vspace{.3cm}

{\bf Remark.} This theorem was first proved, in the
special case $G=SL_n$, by Laumon~\cite{La}. His
argument cannot be generalized to arbitrary semisimple groups.
In the general case, the theorem was proved by Faltings
\cite[theorem~II.5]{Fa}. The proof below seems to be
more elementary than that of Faltings; it is based on nothing but
a few
general results of Symplectic geometry. Another proof of the theorem is
given in [BD]. That proof is more complicated; however, it potentially
leads to a  description of the irreducible components of
$\NN$.
\vspace{.3cm}

\ab We begin with a few general lemmas.
Let $(M_1,\omega_1)$ and $(M_2,\omega_2)$ be complex 
algebraic symplectic
manifolds, and $\pr_i : M_1\times M_2 \to M_i,$ the projections.
We regard $M_1\times M_2$ as a symplectic manifold with symplectic form
$\pr_1^*\omega_1-\pr_2^*\omega_2$, involving the minus-sign on the second
 factor. The following result is a special case of \cite[Prop. 2.7.51]{CG}.
\vspace{2pt}

{\bf Lemma 2.} {\it 
Let $\Lambda_1 \subset M_1$ and
$\Lambda \subset M_1\times M_2$ be smooth algebraic
isotropic subvarieties. Then}
$\pr_2\left(\pr_1^{-1}(\Lambda_1)\cap
\Lambda\right)\sset M_2$ {\it is an isotropic subvariety.
}
\vspace{2mm}

{\sl Proof.}
Set $Y:=\pr_1^{-1}(\Lambda_1)\cap \Lambda$. Simple linear algebra shows that,
 for any $y\in Y$, the image of the tangent
map $(\pr_2)_{_*}: T_yY \to T_{\mbox{{\tiny pr}}_2(y)}M_2$ is isotropic.
We use the characterization of isotropic subvarieties mentioned at
the beginning of the paper.
Let $W\subset \Lambda_2:=\pr_2(Y)$ be an irreducible smooth subvariety.
Observe that the map $\pr_2 :\pr_2^{-1}(W)\cap Y \to W$ is 
{\it surjective}. Hence, there exists a non-empty smooth
Zariski-open dense  subset $Y'\subset 
\left(\pr_2^{-1}(W)\cap Y\right)_{red}$ such
that the restriction $\pr_2: Y' \to W$ has surjective differential
at any point of $Y'$. Therefore the tangent space at the
generic point of $W$ is isotropic. Whence the tangent space at
every point of $W$ is isotropic, by continuity. It follows that
 any smooth subvariety
of $\Lambda_2$ is isotropic, and the lemma follows.
\sq\medskip

\ab Given a manifold $N$ we write $\lambda_N$ for the
canonical
1-form on $T^*N$, usually denoted `$p dq$', such that $d\lambda$ is
the canonical symplectic 2-form on $T^*N$.
Let $f: N_1\to N_2$ be a morphism of smooth algebraic varieties.
Identify $T^*(N_1\times N_2)$ with $T^*N_1\times T^*N_2$ 
via the standard map multiplied by $(-1)$ on the factor
$T^*N_2$.
The canonical
1-form on $T^*(N_1\times N_2)$ becomes
under the above identification  equal to
 $\pr_1^*(\lambda_{N_1}) -\pr_2^*( \lambda_{N_2})$.
We endow $T^*N_1\times T^*N_2$ with the corresponding
symplectic form $\pr_1^*(d\lambda_{N_1}) - \pr_2^*(d\lambda_{N_2})$; it is
induced from the canonical symplectic form 
on $T^*(N_1\times N_2)$ via the  identification.

\ab Introduce the following closed subvariety
$$\small{Y_f= \{(n_1,\alpha_1),(n_2,\alpha_2)\in T^*N_1\times T^*N_2 \;\big|\;
 n_2=f(n_1)\,,\,\alpha_1=0=f^*(\alpha_2)\}\,.}\eqno(3)$$
\vspace{1pt}

{\bf Lemma 4.} {\it The image of $Y_f$ under the second
projection} $\pr_2: T^*N_1\times T^*N_2 \to T^*N_2$
{\it is an isotropic subvariety in} $T^*N_2$.
\vspace{2mm}

{\sl Proof.}
Using the above explained identification
of $T^*(N_1\times N_2)$ with $T^*N_1\times T^*N_2$ involving a sign,
the conormal bundle to the graph of $f$ can be written as the subvariety
\[\Lambda = \{(n_1,\alpha_1),(n_2,\alpha_2)\in T^*N_1\times T^*N_2 \;\;\big|\;\;
 n_2=f(n_1)\,,\,
\alpha_1=f^*(\alpha_2)\}\,.\]
Observe that the canonical 1-form 
$\pr_1^*(\lambda_{N_1}) -\pr_2^*( \lambda_{N_2})$ on
$T^*N_1\times T^*N_2$ vanishes identically on $\Lambda$.
Hence, $\Lambda$ is an
isotropic subvariety, and we may 
apply Lemma  1 
to $M_1= T^*N_1\,,\, M_2=T^*N_2$, and $\Lambda_1=T^*_{N_1}N_1= $
{\it zero-section}, and to $\Lambda$ above.
 Observe now that we have by definition,
$Y_f=\Lambda\cap \pr_1^{-1}(T^*_{N_1}N_1)$.
Hence by Lemma  2  the subvariety $\pr_2(Y_f)$ is isotropic.
\sq

\vspace{2pt}

{\bf Lemma 5.} {\it
If $N_1$ and $N_2$ are smooth algebraic stacks, and
$f: N_1\to N_2$ is a representable morphism of finite type, then the
assertion of Lemma  4  remains valid.}
\vspace{2mm}

{\sl Proof.}
 Due to locality of the claim we may (and will)
assume $N_2$ is quasi-compact.
Let ${\widetilde N_2}$ be a smooth algebraic
variety and  ${\widetilde N_2} \to N_2$ a smooth surjective 
equidimensional morphism.
Set ${\widetilde N_1} := {\widetilde N_2}\times_{N_2}N_1$. 
Note that the set $Y_f \subset T^*N_1\times T^*N_2$ defined in (3)
may be viewed as a
subset in
$T^*N_2\times_{N_2}N_1.$ Therefore
we have 
\[Y_f\times_{N_2}{\widetilde N_2} \,\subset\, 
T^*N_2\times_{N_2}{\widetilde N_2}\times_{N_2}N_1
\,=\, T^*{\widetilde N_2}\times_{N_2}{\widetilde N_1}.\]
We must show that the image of
$Y_f\times_{N_2}{\widetilde N_2}$ is an isotropic subvariety
in $ T^*{\widetilde N_2}$.
Let ${F} : {\widetilde N_1} \to {\widetilde N_2}$ be the natural morphism,
and 
$Y_{_{F}}\subset 
T^*{\widetilde N_2}\times_{_{\widetilde N_2}}{\widetilde
N_1}$
the corresponding subvariety of Lemma  4. Observe
that $T^*{\widetilde N_2}\times_{_{\widetilde N_2}}{\widetilde
N_1}=T^*{\widetilde N_2}\times_{N_2}N_1$ and $Y_f\times_{N_2}{\widetilde N_2}=
Y_F$.
Hence, Lemma  4 applied to $F$
 shows that the image of $Y_{_{F}}\times_{N_2}{\widetilde
N_2}$
in $T^*{\widetilde N_2}$ is isotropic.
The claim follows.
\sq
\vspace{2pt}

\ab Choose a Borel subgroup $B\subset G$ with Lie algebra 
${\mathfrak b}$, and let  ${\mathfrak n}$ denote the nilradical of ${\mathfrak b}$.
Given a field $K$
 of characteristic zero, write $G(K),B(K),{\mathfrak b}(K), $ etc., for
the
corresponding sets of $K$-rational points. The following result seems to
be well-known; it is included here for the reader's convenience.
\vspace{2pt}

{\bf Lemma 6.} {\it For any field $K\supset\C$ and any
$x\in \n(K)$, there exists an element $g\in G(K)$ such that
$\Ad g(x)\in {\mathfrak n}(K)$.}
\vspace{2mm}

{\sl Proof.} Since the Jacobson-Morozov theorem holds for  any field of
 characteristic zero, one may find an
${\mathfrak{s}\mathfrak{l}}_2$-triple
$(x,h,x^-) \subset \g(K)$ associated to the given nilpotent
element $x\in \g(K)$. The eigenspaces of the semisimple endomorphism
$\ad h: \g\to\g$ corresponding to  non-negative eigenvalues span
a parabolic subalgebra $\p_x\subset \g,$ which is
defined over $K$.
Writing $\uu_x$ for the nilradical of $\p_x$, by construction, we have:
$x\in \uu_x(K)$. Clearly, if $\bg_x\subset \p_x$, is
a Borel subalgebra defined over $K$,
and ${\mathfrak n}_x$ is its nilradical, then
$\,x\in \uu_x(K)\subset {\mathfrak n}_x(K)$.

\ab Thus, it suffices to prove that the parabolic
$\p_x$ contains a Borel subalgebra defined over $K$.
To this end, let ${\mathcal{P}}$ denote the partial flag variety
of all parabolics in $\g$ of type $\p_x$. There is a unique
$\p\in {\mathcal{P}}(K)$ such that $\p\supset\bg$, where 
$\bg$ is our fixed Borel subalgebra. Now, the group
$G(K)$ acts transitively on ${\mathcal{P}}(K)$
(this follows easily from the Bruhat decomposition, see [Ja]).
We deduce that there exists $g\in G(K)$ such that
$\Ad g(\p)=\p_x.$ But then,
$\Ad g(\bg)\subset \p_x$ is a Borel subalgebra defined over $K$,
and we are done.\sq
\medskip

\ab  Let $f_1,\ldots,f_r\,,\,(r=\rk\g),$ be a set of
homogeneous free generators of $\C[\g]^G$, the  algebra of
$G$-invariant polynomials on $\g$.  
Let $d_i=\deg f_i$
be the {\it exponents} of $\g$. Following Hitchin \cite{Hi} we put:
$$\hh \;:= \; \Gamma\bigl(X, \Omega_X^{\,\otimes d_1}\bigr)
\;\bigoplus\;\ldots\;\bigoplus\;
\Gamma\bigl(X, \Omega_X^{\,\otimes d_r}\bigr)
\,.$$
This is an affine space of  dimension  equal to $\dim Bun_G$
(at this point it is used, see [Hi], that genus $X$ is $>1$).
Hitchin has defined a morphism
$\pi :T^*Bun_G \to \hh$, by assigning
to any pair $(s, P)\in T^*Bun_G $, where
$P\in Bun_G $ and $s\in T^*_{P}Bun_G
\simeq\Gamma(X,\, {\mathfrak g}_{_{\PP(S)}}\otimes
\Omega^1_{X}),\,$ see (1),
the element
$\pi(s, P)= \oplus_{i=1}^r\; f_i(s)\,\in\,\hh.$
 It is immediate from
the construction
that the global nilpotent variety
is the fiber of $\pi$ over the zero element $0\in\hh$.
\medskip

{\bf Remark.} Hitchin actually worked in the setup of {\it stable} 
Higgs bundles
and
not in the setup of stacks. But his construction of the map $\pi$
extends to the stack setup verbatim. 
We make no  use of any additional properties
of the map $\pi$ established in \cite{Hi}.
\vspace{2pt}

{\bf Lemma 7.} $\;\;\NN$ {\it is an isotropic substack in} $T^*Bun_G$.
\vspace{2mm}

{\sl Proof.} Write $N_2=Bun_B$ for the moduli stack of principal
$B$-bundles. By an old result of Harder [Ha],
any $G$-bundle on a curve has a $B$-reduction, hence
the natural morphism of stacks $f: Bun_B \to Bun_G$ is surjective.

\ab Let $P$ be an  algebraic $G$-bundle  on the curve $X$,
and $s$  a  nilpotent regular section of 
${\mathfrak g}_{_P}\otimes\Omega^1_X$.
Harder showed further, using a key rationality result of 
Steinberg [St], that the bundle $P$
is locally-trivial in the Zariski topology.
Thus, trivializing $P$ on the generic point of $X$, one may 
identify the restriction of $s$ to the generic point
with a nilpotent element of $\g(K)$, where $K=\C(X)$ is the field
of rational functions on $X$.
Hence Lemma 6 implies that
there exists a $B$-reduction 
of  $P$  over the generic point of $X$
 such that $s\in {\mathfrak n}_{_P}\otimes
\Omega^1_X$. Here, a $B$-reduction is a section of the associated bundle
$B\backslash P$. The fibers of $B\backslash P$ being projective
varieties
(isomorphic to $B\backslash G$), any  section of $B\backslash P$ defined
 over the generic point of $X$ extends to the whole of $X$.
Thus, there exists a $B$-reduction  over $X$ 
of  the $G$-bundle $P$    such that $s\in {\mathfrak n}_{_P}\otimes
\Omega^1_X$.

\ab Further, let $k\supset\C$ be a field, and $S=\Spec(k)$.
For any $P\in Bun_B(k)$,  we have:
\[T^*_PBun_B= H^1(X,{\mathfrak b}_{_P})^*= H^0(X,{\mathfrak b}^*_{_P}\otimes
\Omega^1_X)=H^0\left(X,({\mathfrak g}_{_P}/{\mathfrak n}_{_P})\otimes
\Omega^1_X\right).\]
It follows that in the notation of  Lemma  5, for
 $N_1=Bun_B$ and $N_2=Bun_G$,
we have $\NN=\pr_2(Y_f)$. Observe further that 
$Bun_B$ is the union of a countable family of open substacks of
finite type over $Bun_G$ each.
Thus, Lemma  5  implies that $\NN$ is the union
of a {\it countable} family of isotropic substacks.

\ab We claim that if an algebraic stack $S$
is the union
of a  countable  family $\{S_i\}_{i\in{\Bbb N}}$
of  locally closed substacks, then any  quasi-compact substack
of $S$ can be covered by finitely many $S_i$'s.
To prove this, we may assume without loss of generality
that $S$ is itself
quasi-compact. Choose 
 a smooth surjective morphism $\widetilde{S} \onto S$,
where $\widetilde{S}$ is a scheme of finite type.
By our assumptions, there exists a  countable family
$\,\{\widetilde{S}_i\}_{i\in {\mathbb{N}}}\,$ of locally closed
subschemes of $\widetilde{S}$ such that: 
$\widetilde{S}=\cup_i\,\widetilde{S}_i$.
Hence, there exists $n\gg 0$ such that
$\widetilde{S}$ equals the union of the closures
of  $\widetilde{S}_1,\ldots,\widetilde{S}_n$,
thanks to Baer theorem. Hence,
$\widetilde{S}_1\cup \ldots\cup \widetilde{S}_n$
is Zariski dense in $\widetilde{S}$,
and $\,\dim\bigl(\widetilde{S}\,\smallsetminus\,
(\widetilde{S}_1\cup \ldots\cup \widetilde{S}_n)\bigr) < \dim\widetilde{S}.\,$
Arguing by induction
on $\dim \widetilde{S}$ we deduce that 
$\widetilde{S}$ is covered by finitely many
$\widetilde{S}_i$'s. Hence,
for any field $k\supset\C$, the quasi-compact
set $S(k)$ is covered by finitely many
subsets $S_i(k)$.

\ab Thus, we have proved that any open  quasi-compact substack
of $\NN$ can be covered by finitely many isotropic substacks.
This implies that $\NN$ is itself isotropic. \sq
\vspace{2pt}

{\bf Proposition 8.} $\quad\dim\NN=\dim Bun_G$.
\vspace{2mm}

{\sl Proof.}
We observe  that,
 since $Bun_G$ is an equi-dimensional smooth stack,
each irreducible component of $T^*Bun_G$ has dimension $\ge 2\dim
Bun_G$. To see this,
one may replace $Bun_G$ by an  open quasi-compact
substack $Y$, which admits a smooth  surjective morphism 
 $p: \widetilde{Y} \onto Y$,
where $\widetilde{Y}$ is a smooth algebraic variety
and fibers of $p$ are purely $m$-dimensional.
Then,
$T^*Y\times_Y \widetilde{Y}$ is a closed subscheme of
$T^*\widetilde{Y}$ locally defined by $m$ equations. Hence,
for each irreducible component $T^*_j$ of $T^*Y$ we find:
$$\dim T^*_j = \dim(T^*_j\times_Y \widetilde{Y})-m\ge
(\dim T^*\widetilde{Y} -m) -m = 2(\dim \widetilde{Y} - m) =2\dim Y\,.$$
It follows that
each irreducible component of  any fiber of the Hitchin morphism 
$\pi:T^*Bun_G \to \hh$ has dimension $\ge  2\dim Bun_G-\dim\hh=
\dim Bun_G$. But we have proved that
each component of $\NN=\pi^{-1}(0)$ is an isotropic subvariety.
Thus, $\dim\NN=\dim Bun_G$.
\sq
\vskip 2pt

{\bf Remark.} Although the inequality: $\dim T^*Bun_G \ge  2\dim Bun_G$
is no longer true if $X$ has genus one or zero, it has been shown
in [BD] that the stack $\NN$ still has pure 
dimension equal to $\dim Bun_G$.
\medskip

{\bf Corollary 9.} {\it 
Every irreducible component of any fiber of $\pi$ has dimension
equal to}
$\dim Bun_G$. {\it In particular, the Hitchin morphism $\pi$ is flat.}
\vspace{2mm}

{\sl Proof.} There is a natural
$\C^*$-action on $\hh$, such that $t\in\C^*$ acts on the direct
summand $\Gamma(X,\,\Omega_X^{\,\otimes d_i})\subset \hh$ via multiplication
by $t^{d_i}$. 
The  map $\pi: T^*Bun_G \to \hh$
is $\C^*$-equivariant relative to the above-defined
$\C^*$-action on $\hh$ and  the standard $\C^*$-action
on $T^*Bun_G$ by
dilations along the fibers, respectively.
 Clearly, zero is the only fixed point of 
the $\C^*$-action on $\hh$  and, moreover, 
it is contained in the closure of any
other $\C^*$-orbit on $\hh$. It follows, 
since dimension of any irreducible component of any
fiber $\le$ dimension of the
special fiber,
that for any $h\in\hh$ we have:
$\dim\pi\inv(h)\le\dim
\pi\inv(0)$. But for curves of genus $> 1$,
Proposition 8 yields:
$\dim \pi\inv(0) = \dim \NN = \dim Bun_G.\,$
Thus, for any $h\in\hh$ we get:
$\dim\pi\inv(h)\le\dim Bun_G.\,$
On the other hand, the dimension of each irreducible component
of any fiber of the morphism $\pi: T^*Bun_G \to \hh$
is no less than $\dim T^*Bun_G - \dim \hh= \dim Bun_G.\,$
This proves the opposite inequality.
\sq\medskip

\ab The proof of the main theorem is completed by the following stronger
result

\vskip 2pt
{\bf Theorem 10.} {\it The stack $T^*Bun_G$ is a local complete
intersection,
and $\NN$ is a Lagrangian complete intersection in $T^*Bun_G$.}
\vspace{2mm}

{\sl Proof.} The claims being local, we may replace
 $Bun_G$ by an  open quasi-compact
substack $Y$ which admits a presentation
of the form $p: \widetilde{Y} \onto Y$,
where $\widetilde{Y}$ is a smooth algebraic variety
and $p$ is a smooth  surjective morphism 
with fibers of  pure 
dimension $m$. As we have observed in the proof of Proposition 8,
$T^*Y\times_Y \widetilde{Y}$ is a closed subscheme of $T^*\widetilde{Y}$
locally defined by $m$ equations and, moreover, $\dim T^*Y \ge 2\dim Bun_G$.
On the other hand, since fibers of the Hitchin map
$\pi$ have dimension $=\dim Bun_G$,
we find: $\dim T^*Y \leq \dim Bun_G + \dim\hh= 2\dim Bun_G$.
This proves that the stack $T^*Bun_G$ is a local complete
intersection. 

\ab Further,  $\NN$ being the zero fiber of the surjective 
morphism $\pi$, it is defined by $\dim Bun_G$ equations
in $T^*Bun_G$. The equality $\dim\NN =\dim Bun_G$
(Proposition 8) implies that $\NN$ is a complete intersection in
$T^*Bun_G$. Finally, since $\NN$ is equidimensional (Corollary 9),
Lemma 7 implies that $\NN$ is Lagrangian.\sq

\vskip 2mm
\pagebreak[3]
{\footnotesize {\it I am grateful to 
V. Drinfeld for his invaluable help, and also for his 
extreme patience while answering
my numerous foolish questions.}}
\vskip 4pt

\footnotesize{
{\bf V.G.}: University of Chicago, Department of Mathematics,
Chicago, IL
60637, USA;\\ 
\hphantom{x}\ab {\bf ginzburg@math.uchicago.edu}
\end{document}